\documentclass[10pt]{amsart}

\usepackage{tikz}
\usetikzlibrary{positioning}
\usetikzlibrary{arrows}
\usepackage{amssymb,amsmath,latexsym,graphicx}
\usepackage{charter,eucal}
\usepackage{amssymb}
\usepackage{tcolorbox}
\usepackage{amscd}
\usepackage[all,cmtip]{xy}
\usepackage{rotating,multicol,soul}
\usetikzlibrary{decorations.markings}

\setul{}{1.1pt}

\newtheorem{theorem}{Theorem }[section]

\newtheorem{lemma}[theorem]{Lemma}
\newtheorem{observation}[theorem]{Observation}

\newtheorem{remark}[theorem]{Remark}
\newtheorem{corollary}[theorem]{Corollary}
\newtheorem{proposition}[theorem]{Proposition}
\newtheorem{principle}[theorem]{\textsc{Principle}}

\newcommand{\bt}{\begin{theorem}}
\newcommand{\et}{\end{theorem}}
\newcommand{\bmt}{\begin{maintheorem}}
\newcommand{\emt}{\end{maintheorem}}
\newcommand{\bc}{\begin{corollary}}
\newcommand{\bl}{\begin{lemma}}
\newcommand{\ec}{\end{corollary}}
\newcommand{\el}{\end{lemma}}
\newcommand{\bo}{\begin{observation}}
\newcommand{\eo}{\end{observation}}
\newcommand{\bp}{\begin{proposition}}
\newcommand{\ep}{\end{proposition}}
\newcommand{\br}{\begin{remark}}
\newcommand{\er}{\end{remark}}
\newcommand{\bpr}{\begin{principle}}
\newcommand{\epr}{\end{principle}}

\def\I{\mathop{\mathrm{I}}}

\def\I{\mathbf{I}}

\newcommand{\epiarrow}[1]{%
\parbox{#1}{\tikz{\draw[thick,->>](0,0)--(#1,0);}}
}

\def\eop{\hspace*{\fill}$\blacksquare$}

\newcommand{\F}{\mathbb{F}}

\newcommand{\mG}{\mathcal{G}}

\newcommand{\mA}{\mathcal{A}}
\newcommand{\mP}{\mathcal{P}}
\newcommand{\mL}{\mathcal{L}}
\newcommand{\mE}{\mathcal{E}}

\newcommand{\bO}{\mathsf{O}}

\newcommand{\mS}{\mathcal{S}}

\newcommand{\wt}{\widetilde}
\newcommand{\ol}{\overline}

\newcommand{\pp}{\prime\prime}
\newcommand{\ppp}{\prime\prime\prime}
\newcommand{\pppp}{\prime\prime\prime\prime}

\usetikzlibrary{matrix}
\usetikzlibrary{arrows}
\usepackage{pgf}
\usepackage{lscape}
\usepackage{listings}

\usepackage{enumitem}
\title{Epimorphisms of generalized polygons B: The octagons}
\author{Joseph A. Thas and Koen Thas}
\thanks{}
\address{Ghent University, Department of Mathematics, Krijgslaan 281, S25, B-9000 Gent, Belgium}
\email{thas.joseph@gmail.com; koen.thas@gmail.com}
\date{}

\begin{document}

\maketitle

\begin{abstract}
This is the second part of our study of epimorphisms with source a thick generalized $m$-gon and target a thin generalized $m$-gon. We classify the case $m = 8$ when the polygons are finite (in the first part \cite{part3} we handled the cases $m = 3, 4$ and $6$). Then we show that the infinite case is very different, and construct examples which strongly differ from the finite case. A number of general structure theorems are also obtained, and 
we also take a look at the infinite case for general gonality.  
 %We introduce the theory of locally finitely generated generalized polygons and locally finitely chained generalized polygons along the way.\\

\end{abstract}

\begin{tcolorbox}
\tableofcontents
\end{tcolorbox}

\section{Introduction}

In this paper, which is a sequel to \cite{part3}, we will study particular cases of epimorphisms between generalized octagons. Our main motivation of writing the first instance in this series originated from a result of Pasini \cite{P}, which generalizes an older result of Skornjakov \cite{Skorn} and Hughes \cite{H} (the latter two papers handle the case of projective planes). 

\begin{theorem}[Skornjakov--Hughes--Pasini \cite{Skorn,H,P}]
\label{HP}
Let $\alpha$ be a morphism from a thick (possibly infinite) generalized $m$-gon $\mE$ to a thick (possibly infinite) generalized $m$-gon $\mE'$, with $m \geq 3$. If $\alpha$ is surjective, then 
either $\alpha$ is an isomorphism, or each element in $\mE'$ has an infinite fiber in $\mE$.
\end{theorem}
(The mathematical notions will be detailed in the next section.) 
The theorem implies that an epimorphism between finite thick generalized polygons necessarily is an isomorphism, which is quite a surprise, if we would for instance look at this result from the viewpoint of the category of finite groups: the first isomorphism theorem of groups says that for a given epimorphism 
\begin{equation}
\gamma:\ A\ \epiarrow{1cm} \ B,
\end{equation}
we have a natural isomorphism $B \cong A/C$, where $C$ is the kernel of $\gamma$, but finiteness assumptions on $A$ and $B$ by no means imply that $C$ is trivial. (We {\em do} know that the fibers of the elements in $B$ all have size $\vert C \vert$.) So apparently, in the geometric setting of Theorem \ref{HP}, the fact that the source and target polygons are finite puts enough geometric constraints on the morphism to have a ``trivial kernel.''  Van Maldeghem pointed out in \cite[section 4.2.4]{POL} that the thickness assumption is crucial here, and in that same remark he briefly mentions a counter example in the thin case. Wondering what the ``thin version'' of Theorem \ref{HP} would be, was the starting point of \cite{part3}, and in that paper we handled the finite projective planes, finite generalized quadrangles and finite generalized hexagons with thick source and thin target (the latters assumed to come with an order).  

\begin{remark}{\rm 
A local variation on Theorem \ref{HP} by B\"{o}di and Kramer \cite{BK} states that an epimorphism between thick generalized 
$m$-gons ($m \geq 3$) is an isomorphism if and only if its restriction to at least one point row or line pencil is bijective. Later on, Gramlich and Van Maldeghem thoroughly studied epimorphisms from thick generalized $m$-gons to thick generalized $n$-gons with $n < m$ in their works \cite{GramHVM0,GramHVM}, and again, classification results were obtained based on the local nature of the epimorphisms.  }
\end{remark} 

In \cite{part3}, we also considered the infinite case, and showed that the finite formulation does not naturally generalize to the infinite case (through the construction of counter examples). For that matter, we introduced locally finitely chained and locally finitely generated polygons. 

In this paper, we handle the case of generalized octagons (which is more tedious than the other cases), and we also consider the infinite case by constructing (counter) examples in the spirit of \cite{part3}.  

%A second motivation came from our recent paper \cite{part2}, which is a second instalment in a series of papers \cite{part1,part2} about covers of generalized quadrangles. In that study, a question naturally popped up involving morphisms between thick finite generalized quadrangles which surjectively map some geometrical hyperplane of the source quadrangle to a geometrical hyperplane of the target quadrangle. The precise statement is described in the second last section of this paper.  

\subsection*{Synopsis of the present paper}

In section \ref{basic}, we introduce some notions which will be frequently used throughout the paper. In section \ref{syn}, we remind the reader of the main results of \cite{part3} regarding  epimorphisms from thick finite generalized $m$-gons with $m \in \{ 3, 4, 6 \}$ to thin generalized $m$-gons with an order. In the subsequent section \ref{main}, we formulate and prove the missing octagonal case (Theorem \ref{JATGO}), which comes with a byproduct for locally finite octagons.   
Then, in section \ref{chain} we  
provide infinite examples of (known) generalized octagons which yield classes of epimorphisms which agree with  our results in the finite case. This is done through the theory of  locally finitely generated generalized polygons and locally finitely chained generalized polygons, which was introduced in the predecessor of the present paper. 
In section \ref{counter} however, we freely construct classes of examples of epimorphisms which show that the infinite case necessarily behaves differently (for all finite gonalities $n \geq 3$). The point is that there are no restrictions on the order $(s',1)$ of the thin target polygon. 
%Finally, in section \ref{lf} we explain a corollary of the proof of Theorem \ref{JATGO} for locally finite octagons. 

\section{Some basic definitions}
\label{basic}

We summarize a number of definitions which we will need in due course. 

\subsection{Generalized polygons}

Let $\Gamma = (\mP,\mL,\I)$ be a point-line geometry, and let $m$ be a positive integer at least $2$. We say that $\Gamma$ is a {\em weak generalized $m$-gon} if any two elements in $\mP \cup \mL$ are contained in at least one ordinary sub $m$-gon (as a subgeometry of $\Gamma$), and if $\Gamma$ does not contain ordinary sub $k$-gons with $2 \leq k < m$. For $m = 2$ every point is incident with every line. 

If $m \geq 3$, we say $\Gamma$ is a {\em generalized $m$-gon} if furthermore $\Gamma$ contains an ordinary sub $(m + 1)$-gon as a subgeometry. Equivalently, a weak generalized $m$-gon with $m \geq 3$ is a generalized $m$-gon if it is {\em thick}, meaning that every point is incident with at least three distinct lines and every line is incident with at least three distinct points. A weak generalized $m$-gon is {\em thin} if it is not thick; in that case, we also speak of {\em thin generalized $m$-gons}. If we do not specify $m$ (the ``gonality''), we speak of {\em (weak) generalized polygons}. Note that thick generalized $2$-gons (or {\em generalized digons}) do not contain ordinary $3$-gons as a subgeometry.  
 
It can be shown that generalized polygons have an {\em order} $(u,v)$: there exists positive integers $u \geq 2$ and $v \geq 2$ such that each 
point is incident with $v + 1$ lines and each line is incident with $u + 1$ points. We say that a weak generalized polygon is {\em finite} if its number of points and lines is finite | otherwise it is {\em infinite}.  If a thin weak generalized polygon has an order $(1,u)$ or $(u,1)$ it is called a {\em thin 
generalized polygon} of order $(1,u)$ or $(u,1)$. 

Note that the generalized $3$-gons are precisely the (axiomatic) projective planes. Generalized $4$-gons, resp. $6$-gons, resp. $8$-gons are also called {\em generalized quadrangles}, resp. {\em hexagons}, resp. {\em octagons}.

\subsection{Sub polygons}

We say that $\Gamma' = (\mP',\mL',\I')$ is a {\em sub generalized $m$-gon} of the generalized $m$-gon $\Gamma$, $m \geq 3$, if $\mP' \subseteq \mP$, $\mL' \subseteq \mL$, and if $\I'$ is the induced incidence coming from $\I$.

\subsection{Morphisms and epimorphisms}

A {\em morphism} from a weak generalized polygon $\Gamma = (\mP,\mL,\I)$ to a weak generalized polygon $\Gamma = (\mP',\mL',\I')$ is a map $\alpha: \mP \cup \mL \mapsto \mP' \cup \mL'$ which maps points to points, lines to lines and which preserves the incidence relation (note that we do not ask the gonalities to be the same). We say that a morphism $\alpha$ is an {\em epimorphism} if $\alpha(\mP) = \mP'$ and $\alpha(\mL) = \mL'$. Contrary to Gramlich and Van Maldeghem \cite{GramHVM0,GramHVM}, we do not ask surjectivity onto the set of flags of $\Gamma'$ (the incident point-line pairs of $\Gamma'$). 

Note that in categorical language, an {\em epimorphism} is any morphism which is right-cancellative. 
In the category of sets, this is trivially equivalent to asking that the morphism (map) is surjective. Since morphisms between generalized polygons are defined by the underlying maps between the point sets and line sets, it follows that in the categorical sense, epimorphisms between polygons are indeed as above.  

If an epimorphism is injective, and if the inverse map is also a morphism, then we call it an {\em isomorphism}.

\subsection{Doubling}

Let $\Gamma = (\mP, \mL, \I)$ be a (not necessarily finite) generalized $n$-gon of order $(s,s)$ (for $n = 3$, projective planes of order $(1,1)$ are allowed). Define the {\em double of $\Gamma$} as the generalized $2n$-gon $\Gamma^\Delta$ which arises by letting its point set be $\mP \cup \mL$, and letting its line set be the flag set of $\Gamma$ (the set of incident point-line pairs).  Its parameters are $(1, s)$. The full automorphism group of $\Gamma^\Delta$ is isomorphic to the  group consisting of all automorphisms and dualities (anti-automorphisms) of $\Gamma$. Sometimes we prefer to work in the point-line dual of $\Gamma^\Delta$, but we use the same notation (while making it clear in what setting we work). This is what we will do in this section.
Vice versa, if $\Gamma'$ is a thin generalized $2n$-gon of order $(1,s)$, then it is isomorphic to the double $\Gamma^\Delta$ of a generalized $n$-gon $\Gamma$ of order $(s,s)$. \\

\section{Synopsis of known results}
\label{syn}

In this section we summarize some of the results of \cite{part3}. Specifically, we list the theorems obtained on epimorphisms from finite projective planes, finite generalized quadrangles and finite generalized hexagons to thin planes, quadrangles and hexagons, respectively. In the general framework of classifying epimorphisms with source a finite thick generalized $n$-gon and target a thin $n$-gon (with an order), only the octagons (case $n = 8$) remain, due to a result of Feit and Higman \cite{FEHI}.

\begin{theorem}[The planes \cite{part3}]
\label{GT}
Let $\Phi$ be an epimorphism of a thick projective plane $\mP$ onto a thin projective plane $\Delta$ of order $(1,1)$.  Then exactly two classes of epimorphisms $\Phi$ occur (up to a suitable permutation of the points of $\Delta$), and they are described as follows. 

\begin{itemize}
\item[{\rm(a)}] 
The points of $\Delta$ are $\overline{a}, \overline{b}, \overline{c}$, with $\overline{a} \sim \overline{b} \sim \overline{c} \sim \overline{a}$, and put $\Phi^{-1} (\overline{x}) = \wt X$, with $\overline{x} \in \lbrace \overline{a}, \overline{b}, \overline{c} \rbrace$.

   Let $(\wt{A}, \wt{B})$, with $\wt{A} \ne \emptyset \ne \wt{B}$, be a partition of the set of all points incident with a line $L$ of $\mP$. Let  $\wt{C}$ consist of the points not incident with $L$. Furthermore, $\Phi^{-1}(\overline{a}\overline{b}) = \{ L\}$, $\Phi^{-1}(\overline{b}\overline{c})$ is the set of all lines distinct from $L$ but incident with a point of $\wt{B}$ and $\Phi^{-1}(\overline{a}\overline{c})$ is the set of all lines distinct from $L$ but incident with a point of $\wt{A}$.
   \item[{\rm(b)}] 
   The dual of (a).
\end{itemize}
\end{theorem}

\begin{theorem}[The quadrangles \cite{part3}]
\label{JATGQ}
Let $\Phi$ be an epimorphism of a thick generalized quadrangle $\mathcal{S}$ of order $(s, t)$ onto a grid $\mathcal{G}$. Let $\mathcal{G}$ have order $(s^\prime, 1)$. Then $s^\prime = 1$ and exactly two classes of epimorphisms $\Phi$ occur (up to a suitable permutation of the points of $\mathcal{G})$.

\begin{itemize}
\item[{\rm(a)}] The points of $\mathcal{G}$ are $\overline{a}, \overline{b}, \overline{c}, \overline{d}$, with $\overline{a} \sim \overline{b} \sim \overline{c} \sim \overline{d} \sim \overline{a} $, and put $\Phi^{-1} (\overline{x}) = \wt X$, with $\overline{x} \in \lbrace \overline{a}, \overline{b}, \overline{c}, \overline{d} \rbrace$.

   Let $(\wt{A}, \wt{B})$, with $1 \le \vert \wt A \vert \le s, 1 \le \vert \wt B \vert \le s$, be a partition of the set of all points incident with a line $L$ of $\mathcal{S}$. Let  $\wt{C}$ consist of the points not incident with $L$ but collinear with a point of $\wt{B}$, and let $\wt{D}$ consist of the points not incident with $L$ but collinear with a point of $\wt{A}$. Further, $\Phi^{-1}(\overline{a}\overline{b}) = \{ L \}$, $\Phi^{-1}(\overline{b}\overline{c})$ is the set of all lines distinct from $L$ but incident with a point of $\wt{B}$, $\Phi^{-1}(\overline{a}\overline{d})$ is the set of all lines distinct from $L$ but incident with a point of $\wt{A}$ and $\Phi^{-1}(\overline{c}\overline{d})$ consists of all lines incident with at least one point of $\wt{C}$ and at least one point of $\wt{D}$.
\item[{\rm(b)}] The dual of (a).
\end{itemize}
\end{theorem}

  \begin{theorem}[The hexagons \cite{part3}]
   \label{JATGH}
   Let $\Phi$ be an epimorphism of a thick generalized hexagon $\mS$ of order $(s, t)$ onto a thin generalized hexagon $\mG$ of order $(s^\prime, 1)$. Then $s^\prime = 1$ and exactly two classes of epimorphisms $\Phi$ occur (up to a suitable permutation of the points of $\mG$).
   
   \begin{itemize}
   \item[{\rm (a)}] The points of $\mG$ are $\ol{a}, \ol{b}, \ol{c}, \ol{d}, \ol{e}, \ol{f}$, with $\ol{a} \sim \ol{b} \sim \ol{c} \sim \ol{d} \sim \ol{e} \sim \ol{f} \sim \ol{a}$, and put $\Phi^{-1}(\ol{x}) = \wt{X}$, with $ \ol{x} \in \lbrace \ol{a}, \ol{b}, \ol{c}, \ol{d}, \ol{e}, \ol{f} \rbrace$.
   
   Let $(\wt{C}, \wt{B}), 1 \le \vert \wt{C} \vert \le s, 1 \le \vert \wt{B} \vert \le s$, be a partition of the set of all points incident with some line $L$ of $\mS$. Let $\wt{D}$ consist of the points not incident with $L$ but collinear with a point of $\wt{C}$, let $\wt{A}$ consist of the points not incident with $L$ but collinear with a point of $\wt{B}$, let $\wt{E}$ consist of the points not in $\wt{C} \cup \wt{D}$ but collinear with a point of $\wt{D}$, and let $\wt{F}$ consist of the points not in $\wt{A} \cup \wt{B}$ but collinear with a point of $\wt{A}$. Further, $\Phi^{-1}(\ol{b}\ol{c}) = \lbrace L \rbrace$, $\Phi^{-1}(\ol{c}\ol{d})$ is the set of all lines distinct from $L$ but incident with a point of $\wt{C}$, $\Phi^{-1}(\ol{a}\ol{b})$ is the set of all lines distinct from $L$ but incident with a point of $\wt{B}$, $\Phi^{-1}(\ol{d}\ol{e})$ is the set of all lines distinct from the lines of $\Phi^{-1}(\ol{d}\ol{c})$ but incident with a point of $\wt{D}$, $\Phi^{-1}(\ol{f}\ol{a})$ is the set of all lines distinct from the lines of $\Phi^{-1}(\ol{a}\ol{b})$ but incident with a point of $\wt{A}$, $\Phi^{-1}(\ol{d}\ol{e})$ is the set of all lines distinct from the lines 
   of $\Phi^{-1}(\ol{c}\ol{d})$ but incident with a point of $\widetilde{D}$,
   and $\Phi^{-1}(\ol{f}\ol{e})$ is the set of all lines not in $\Phi^{-1}(\ol{f}\ol{a})$ but incident with a point of $\wt{F}$ (that is, the set of all lines not in $\Phi^{-1}(\ol{e}\ol{d})$ but incident with a point of $\wt{E}).$
   
   \item[{\rm (b)}] The dual of (a).
  
  \end{itemize}
  \end{theorem}

In the next section, we handle the last remaining (and most tedious) case: the octagons.

\section{Epimorphisms to thin octagons}
\label{main}

\begin{theorem}
\label{JATGO}
Let $\Phi$ be an epimorphism of a thick generalized octagon $\mS$ of order $(s, t)$ onto a thin generalized octagon $\mG$ of order $(s^\prime, 1)$. Then $s^\prime = 1$ and exactly two classes of epimorphisms $\Phi$ occur (up to a suitable permutation of the points of $\mG$).

\begin{itemize}
\item[{\rm (a)}] The points of $\mG$ are $\ol{a}, \ol{b}, \ol{c}, \ol{d}, \ol{e}, \ol{f}, \ol{g}, \ol{h}$, with $\ol{a} \sim \ol{b} \sim \ol{c} \sim \ol{d} \sim \ol{e} \sim \ol{f} \sim \ol{g} \sim \ol{h} \sim \ol{a}$, and put $\Phi^{-1}(\ol{x}) = \wt{X}$, with $\ol{x} \in \lbrace \ol{a}, \ol{b}, \ol{c}, \ol{d}, \ol{e}, \ol{f}, \ol{g}, \ol{h} \rbrace$.

Let $(\wt{C}, \wt{B}), 1 \le \vert \wt{C} \vert \le s, 1 \le \vert \wt{B} \vert \le s$, be a partition of the set of all points incident with a line $L$ of $\mS$. Let $\wt{D}$ consist of the points not incident with $L$ but collinear with a point of $\wt{C}$, let $\wt{A}$ consist of the points not incident with $L$ but collinear with a point of $\wt{B}$, let $\wt{E}$ consist of the points not in $\wt{C} \cup \wt{D}$ but collinear with a point of $\wt{D}$, let $\wt{H}$ consist of the points not in $\wt{A} \cup \wt{B}$ but collinear with a point of $\wt{A}$, let $\wt{F}$ consist of the points not in $\wt{D} \cup \wt{E}$ but collinear with a point of $\wt{E}$, and let $\wt{G}$ consist of the points not in $\wt{A} \cup \wt{H}$ but collinear with a point of $\wt{H}$.

Further, $\Phi^{-1}(\ol{c}\ol{b}) = \lbrace L \rbrace$, $\Phi^{-1}(\ol{c}\ol{d})$ is the set of all lines distinct from $L$ but incident with a point of $\wt{C}$, $\Phi^{-1}(\ol{a}\ol{b})$ is the set of all lines distinct from $L$ but incident with a point of $\wt{B}$, $\Phi^{-1}(\ol{d}\ol{e})$ is the set of all lines distinct from the lines of $\Phi^{-1}(\ol{d}\ol{c})$ but incident with a point of $\wt{D}$, $\Phi^{-1}(\ol{a}\ol{h})$ is the set of all lines distinct from the lines of $\Phi^{-1}(\ol{a}\ol{b})$ but incident with a point of $\wt{A}$, $\Phi^{-1}(\ol{e}\ol{f})$ is the set of all lines distinct from the lines of $\Phi^{-1}(\ol{d}\ol{e})$ but incident with a point of $\wt{E}$, $\Phi^{-1}(\ol{g}\ol{h})$ is the set of all lines distinct from the lines of $\Phi^{-1}(\ol{a}\ol{h})$ but incident with a point of $\wt{H}$, and $\Phi^{-1}(\ol{g}\ol{f})$ is the set of lines not in $\Phi^{-1}(\ol{h}\ol{g})$ but incident with a point of $\wt{G}$ (that is, the set of all lines not in $\Phi^{-1}(\ol{e}\ol{f})$ but incident with a point of $\wt{F}$).

\item[{\rm (b)}] The dual of (a).
\end{itemize}
\end{theorem}

{\em Proof}. \quad We proceed in a number of steps. We first explain a connection between thin generalized octagons and generalized quadrangles.
  
Let $\mG$ be a thin generalized octagon of order $(s^\prime, 1), s^\prime \ge 1$. Define an equivalence relation $\sim$ on the set of lines of $\mG$: \ul{$L \sim M$ if $\widetilde{\bf{d}}$$(L, M)$ = even,} where $\widetilde{\bf{d}}(\cdot,\cdot)$ is the distance in the line graph of $\mG$.  
Let $U, V$ be the equivalence classes. Call the elements of $U$ {\it points}, the elements of $V$ {\it lines}, and $L \in U$ is {\it incident} with $M \in V$ if and only if $\widetilde{\bf{d}}$$(L, M) = 1$. Then this incidence structure is a generalized quadrangle $\widetilde{\mS}$ of order $s^\prime$. Let $x$ be a point of $\mG$, and let $L \in U$ and $M \in V$ be the lines of $\mG$ which are incident with $x$. Then the point $x$ can be identified with the flag $(L, M)$ of the quadrangle $\widetilde{\mS}$.\\

Below, if we use the notation $\mathbf{d}(\cdot,\cdot)$, we express the distance measured in the incidence graph.\\

\quad {\rm (1)} {\bf Let $L$ be any line of $\mS$. Then $\Phi$ maps the set of all points incident with $L$ onto the set of all points incident with $\Phi(L)$}.

{\em Proof of} (1).\quad 
Let $L$ be any line of $\mS$, let $\Phi(L) = \ol{L}$, let $\widehat{L}$ be the set of all points incident with $L$, and let $\widehat {\ol{L}}$ be the set of all points incident with $\ol{L}$. Assume, by way of contradiction, that $\Phi(\widehat{L}) \not = \widehat{\ol{L}}$. Let $\ol{x} \: \mathbf{I} \ol{L}, \ol{x} \not \in \Phi(\widehat{L})$. Let {\bf{d}}$(\ol{x}, \ol{y})$ = 6, {\bf{d}}$(\ol{y}, \ol{L})$ = 7 and $y \in \Phi^{-1}(\ol{y})$. If $z \: \mathbf{I} \: L$ such that {\bf{d}}$(z, y) \le 6$, then {\bf{d}}$(z, y)$ = 6, {\bf{d}}$(y, L)$ = 7 and $\Phi(z) = \ol{x}$. So $\ol{x} \in \Phi(\widehat{L})$, a contradiction.\\

\quad {\rm (2)} {\bf Let $x$ be any point of $\mS$, and let $\Phi(x) = \ol{x}$. Then the image of the set of lines incident with $x$, is the set of lines incident with $\ol{x}$}.

{\em Proof of} (2).\quad 
Dualize proof of (1).\\

\quad {\rm (3)} {\bf $s^\prime = 1$}.

{\em Proof of} (3).\quad
Let $A, B$ be the two systems of lines of $\mG$, let $\Phi^{-1}(A)$ be the union of all inverse images of the elements of $A$, and let $\Phi^{-1}(B)$ be the union of all inverse images of the elements of $B$. If $x$ is a point of $\mS$, then $A_x$ is the set of all lines of $\Phi^{-1}(A)$ incident with $x$ and $B_x$ is the set of all lines of $\Phi^{-1}(B)$ incident with $x$.

Let $c$ be a point of $\mS$ and let $M$ be a line incident with $c$. Let $\ol{c} = \Phi(c)$, let $\ol{d}$ be a point at distance 8 from $\ol{c}$, and let $\Phi(d) = \ol{d}$. Then {\bf{d}}$(c, d)$ = 8. The line $\ol{M} = \Phi(M)$ is incident with $\ol{c}$; say $\ol{M} \in B$. Let $\ol{d} \: \mathbf{I} \: \ol{U}$, {\bf{d}}$(\ol{U}, \ol{M})$ = 6 and $d \: \mathbf{I} \: U$, {\bf{d}}$(U, M)$ = 6. Then $U \in \Phi^{-1}(\ol{U})$ and $\ol{U} \in A$. It easily follows that $\vert B_c \vert = \vert A_d \vert$. Similarly $\vert A_c \vert = \vert B_d \vert$.

From now on, assume, by way of contradiction that $s^\prime > 1$.

Let $(f, F)$ and $(g, G)$ be the point-line flags of the generalized quadrangle $\wt \mS$ which correspond to $\ol{c}$ and $\ol{d}$. Then {\bf{d}}$(\ol{c}, \ol{d})$ = 8 is equivalent with $f \not \sim g$ and $F \not \sim G$ in $\wt \mS$. Let $(h, H)$ be a flag of $\wt \mS$ with $f \not \sim h \not \sim g$ and $F \not \sim H \not \sim G$ (such a flag exists as $s^\prime > 1$). If $\ol{r}$ is the point of $\mG$ which corresponds to $(h, H)$, then {\bf{d}}$(\ol{c}, \ol{r})$ = {\bf{d}}$(\ol{d}, \ol{r})$ = 8. If $r \in \Phi^{-1}(\ol{r})$, then {\bf{d}}$(c, r)$ = {\bf{d}}$(d, r)$ = 8. So $\vert A_c \vert = \vert B_r \vert = \vert A_d \vert$ and similarly $\vert B_c \vert = \vert B_d \vert$. Hence $\vert A_c \vert = \vert A_d \vert = \vert B_c \vert = \vert B_d \vert$.

Let $m$ and $n$ be points at distance 8 in $\mS$. If {\bf{d}}$(\ol{m}, \ol{n})$ = 8, with $\Phi(m) = \ol{m}$ and $\Phi(n) = \ol{n}$, then we know that $\vert B_m \vert = \vert A_m \vert =  \vert B_n \vert = \vert A_n \vert$. 

Now let {\bf{d}}$(\ol{m}, \ol{n})$ = 6, and let $\ol{m} \: \mathbf{I} \: \ol{M}, \ol{n} \: \mathbf{I} \: \ol{N}$, {\bf{d}}$(\ol{M}, \ol{N})$ = 4, {\bf{d}}$(\ol{m}, \ol{u})$ = {\bf{d}}$(\ol{u}, \ol{v})$ = {\bf{d}}$(\ol{v}, \ol{n})$ = 2. Let $\ol{l} \: \mathbf{I} \: \ol{u}\ol{v}, \ol{l} \not = \ol{u}, \ol{l} \not = \ol{v}$. Further, let $\ol{r}$ be such that {\bf{d}}$(\ol{l}, \ol{r})$ = 4, {\bf{d}}$(\ol{u}, \ol{r})$ = {\bf{d}}$(\ol{v}, \ol{r})$ = 6. Then {\bf{d}}$(\ol{m}, \ol{r})$ = {\bf{d}}$(\ol{n}, \ol{r})$ = 8. So $\vert A_m \vert = \vert A_r \vert = \vert B_r \vert = \vert B_m \vert$ and $\vert A_n \vert = \vert A_r \vert = \vert B_r \vert = \vert B_n \vert$. Hence $\vert A_m \vert = \vert B_m \vert = \vert A_n \vert = \vert B_n \vert$.

Next, let {\bf{d}}$(\ol{m}, \ol{n})$ = 4. Let $(f, F)$ and $(g, G)$ be the point-line flags of the generalized quadrangle $\wt \mS$ which correspond to $\ol{m}$ and $\ol{n}$, where $g \: \mathbf{I} \: F$ (in $\wt \mS$). Let $(h, H)$ be a flag in $\wt \mS$ with $h \not \sim f, h \not \sim g, H \not \sim F, H \not \sim G$ (as $s^\prime > 1$ this flag exists). If $\ol{r}$ is the point of $\mG$ which corresponds to $(h, H)$, then {\bf{d}}$(\ol{m}, \ol{r})$ = {\bf{d}}$(\ol{n}, \ol{r})$ = 8. Let $r \in \Phi^{-1}(\ol{r})$; then $\vert A_m \vert = \vert A_r \vert = \vert B_r \vert = \vert B_m \vert$ and $\vert A_n \vert = \vert A_r \vert = \vert B_r \vert = \vert B_n \vert$, so $\vert A_m \vert = \vert B_m \vert = \vert B_n \vert = \vert A_n \vert$.

Finally, let {\bf{d}}$(\ol{m}, \ol{n})$ = 2. Choose $\ol{r}$ such that {\bf{d}}$(\ol{m}, \ol{r})$ = {\bf{d}}$(\ol{n}, \ol{r})$ = 8, {\bf{d}}$(\ol{r}, \ol{m}\ol{n})$ = 7. As before it easily follows that $\vert A_m \vert = \vert A_n \vert = \vert B_m \vert = \vert B_n \vert$.

Hence for $m$ and $n$ points of $\mS$ at  distance 8, we always have $\vert A_m \vert = \vert A_n \vert = \vert B_m \vert = \vert B_n \vert$.

Next, consider points $m$ and $n$ of $\mS$ at distance 6, 4 or 2. Clearly there is a point $r$ with {\bf{d}}$(m, r)$ = {\bf{d}}$(n, r)$ = 8. Again it easlly follows that $\vert A_m \vert = \vert A_n \vert = \vert B_m \vert = \vert B_n \vert$.

We conclude that for any two points $m$ and $n$ of $\mS$ we have 
\begin{equation}
\vert A_m \vert = \vert A_n \vert = \vert B_m \vert = \vert B_n \vert.
\end{equation}

Let $r$ be a point of $\mS$. Let $r \: \mathbf{I} \: L_1, r \: \mathbf{I} \: L_2, L_1 \not = L_2$ and $L_1, L_2 \in A_r$ (remark that $\vert A_r \vert = \frac{t + 1}{2}$). Let $\ol{r} = \Phi(r), \Phi(L_1) = \Phi(L_2) = \ol{L}, \ol{u}_1 \: \mathbf{I} \: \ol{L}, \ol{u}_1 \not = \ol{r}$, and $u_1 \: \mathbf{I} \: L_1$ with $\ol{u}_1 = \Phi(u_1)$. Further, let $\ol{u}_2 \: \mathbf{I} \: \ol{L}, \ol{u}_2 \not = \ol{u}_1$, and $u_2 \: \mathbf{I} \: L_2$ with $\ol{u}_2 = \Phi(u_2)$. Next, let $N_1 \in B_{u_2}, \ol{N}_1 = \Phi(N_1), \ol{v}_1 \: \mathbf{I} \: \ol{N}_1, \ol{v}_1 \not = \ol{u}_2$, and $v_1 \: \mathbf{I} \: N_1$ with $\ol{v}_1 = \Phi(v_1)$. Finally, let $N_2 \in A_{v_1}, \ol{N}_2 = \Phi(N_2), \ol{v}_2 \: \mathbf{I} \: \ol{N}_2, \ol{v}_2 \not = \ol{v}_1$, and $v_2 \: \mathbf{I} \: N_2$ with $\ol{v}_2 = \Phi(v_2)$.

Let $R_1 \in B_{u_1}$ and let {\bf{d}}$(R_1, R_2)$ = 6 with $v_2 \: \mathbf{I} \: R_2$. Hence $R_1$ determines $R_2$. Then $\Phi(R_2) = \ol{R_2} = \ol{N}_2$. Hence $R_2 \in A_{v_2}$. It follows that $\vert B_{u_1} \vert \le \vert A_{v_2}  \setminus \lbrace N_2 \rbrace \vert$. So $\vert B_{u_1} \vert < \vert A_{v_2} \vert$, clearly a contradiction.\\

\quad {\rm (4)} {\bf Let $s^\prime = 1$. Then not for any given points $\ol{x}$ and $\ol{y}$ in $\mG$, with {\bf{d}}$(\ol{x}, \ol{y})$ = 6, there exist points $x \in \Phi^{-1}(\ol{x})$ and $y \in \Phi^{-1}(\ol{y})$ with {\bf{d}}$(x, y)$ = 8}.

{\bf{Remark}.} {\rm Let {\bf{d}}$(\ol{u}, \ol{v})$ = 8, with $\ol{u} = \Phi(u)$ and $\ol{v} = \Phi(v)$ points of $\mG$. Then by (3) $\vert B_u \vert = \vert A_v \vert$ and $\vert A_u \vert = \vert B_v \vert$}.

{\em Proof of} (4).\quad
Assume, by way of contradiction, that for any given points $\ol{x}$ and $\ol{y}$ in $\mG$, with {\bf{d}}$(\ol{x}, \ol{y})$ = 6, there exist points $x \in \Phi^{-1}(\ol{x})$ and $y \in \Phi^{-1}(\ol{y})$ with {\bf{d}}$(x, y)$ = 8.

Let $\ol{a}, \ol{b}, \ol{c}, \ol{d}, \ol{e}, \ol{f}, \ol{g}, \ol{h}$ be the points of $\mG$, with $\ol{a} \sim \ol{b} \sim \ol{c} \sim \ol{d} \sim \ol{e} \sim \ol{f} \sim \ol{g} \sim \ol{h} \sim \ol{a}$. Let $\ol{a} = \Phi(a)$. Choose $h, g, f, e, N$ such that $\ol{h} = \Phi(h), \ol{g} = \Phi(g), \ol{f} = \Phi(f), \ol{e} = \Phi(e), \ol{d}\ol{e} = \Phi(N)$ and $e \: \mathbf{I} \: N$, with $a \sim h \sim  g \sim f \sim e$. Let $d \: \mathbf{I} \: N$ with {\bf{d}}$(a, d)$ = 6, and let {\bf{d}}$(a, b)$ = {\bf{d}}$(b, c)$ = {\bf{d}}$(c, d)$ = 2. Then $\ol{b} = \Phi(b), \ol{c} = \Phi(c), \ol{d} = \Phi(d)$. Let $A = \lbrace \ol{b}\ol{c}, \ol{d}\ol{e}, \ol{f}\ol{g}, \ol{h}\ol{a} \rbrace$ and $B = \lbrace \ol{a}\ol{b}, \ol{c}\ol{d}, \ol{e}\ol{f}, \ol{g}\ol{h} \rbrace$. We have $\vert A_a \vert = \vert B_e \vert$ and $ \vert B_a \vert = \vert A_e \vert, \vert A_b \vert = \vert B_f \vert$ and $\vert B_b \vert = \vert A_f \vert, \vert A_c \vert = \vert B_g \vert$ and $\vert B_c \vert = \vert A_g \vert, \vert A_d \vert = \vert B_h \vert$ and $\vert B_d \vert = \vert A_h \vert$.

Let $\Phi(a^\prime) = \ol{a}$. As {\bf{d}}$(\ol{a}, \ol{e})$ = 8, we have {\bf{d}}$(a^\prime, e)$ = 8. Hence $\vert B_a \vert = \vert A_e \vert = \vert B_{a^\prime} \vert$, similarly $\vert A_a \vert = \vert A_{a^\prime} \vert$.

By assumption there are points $a^\prime$ and $d^\prime$, with $\Phi(a^\prime) = \ol{a}, \Phi(d^\prime) = \ol{d}$ and {\bf{d}}$(a^\prime, d^\prime)$ = 8. Let $d^\prime \sim l \sim n \sim m \sim a^\prime$ and assume that $\Phi(ma^\prime) \in A$ (points can be chosen in such a way). Then $\Phi(m) = \ol{a}, \Phi(n) = \ol{b}$ and $\Phi(l) = \ol{c}$. So $\Phi(ld^\prime) = \ol{c}\ol{d} \in B$.  By letting $m$ vary such that $\Phi(ma') \in A$, we see that 
$\vert B_{d^\prime} \vert \ge \vert A_{a^\prime} \vert$, so $\vert B_d \vert \ge \vert A_a \vert$. Also $\vert B_d \vert = \vert A_h \vert, \vert A_a \vert = \vert B_e \vert$ and so $\vert B_d \vert \ge \vert B_e \vert $.

Similarly one obtains $\vert B_e \vert \ge \vert B_d \vert$, and so $\vert B_e \vert = \vert B_d \vert$.

Consequently $\vert B_a \vert = \vert B_b \vert = \vert B_c \vert = \vert B_d \vert = \vert B_e \vert = \vert B_f \vert = \vert B_g \vert = \vert B_h \vert$, and similarly $\vert A_a \vert = \vert A_b \vert = \vert A_c \vert = \vert A_d \vert = \vert A_e \vert = \vert A_f \vert = \vert A_g \vert = \vert A_h \vert$. As $\vert B_e \vert = \vert A_a \vert$, it now follows that $\vert A_u \vert = \vert B_u \vert = \frac {(t + 1)}{2}$ for each point $u$ of $\mS$.

Such as in the last part of (3) we now find a contradiction.\\

\quad {\rm (5)} {\bf For $s^\prime = 1$ two (mutually dual) cases occur}.

{\em Proof of} (5).\quad
By (4) we know that not for any given points $\ol{x}$ and $\ol{y}$ in $\mG$ with {\bf{d}}$(\ol{x}, \ol{y})$ = 6, there exist points $x \in \Phi^{-1}(\ol{x})$ and $y \in \Phi^{-1}(\ol{y})$ with {\bf{d}}$(x, y)$ = 8.

Let $\ol{a}, \ol{b}, \ol{c}, \ol{d}, \ol{e}, \ol{f}, \ol{g}, \ol{h}$ be the points of $\mG$, with $\ol{a} \sim \ol{b} \sim \ol{c} \sim \ol{d} \sim \ol{e} \sim \ol{f} \sim \ol{g} \sim \ol{h} \sim \ol{a}$, let $A = \lbrace \ol{b}\ol{c}, \ol{d}\ol{e}, \ol{f}\ol{g}, \ol{h}\ol{a} \rbrace$ and $B = \lbrace \ol{a}\ol{b}, \ol{c}\ol{d}, \ol{e}\ol{f}, \ol{g}\ol{h} \rbrace$, and let $\Phi^{-1}(\ol{x}) = \wt{X}$ for $\ol{x} \in \lbrace \ol{a}, \ol{b}, \ol{c}, \ol{d}, \ol{e}, \ol{f}, \ol{g}, \ol{h} \rbrace$.

Without loss of generality, we assume that there do not exist points $a^\prime \in \wt{A}$ and $d^\prime \in \wt{D}$ for which {\bf{d}}$(a^\prime, d^\prime)$ = 8 (and so {\bf{d}}$(a^\prime, d^\prime)$ = 6).\\

\fbox{
\quad{\bf CASE (a)} {\it Assume that $b_1 \sim b_2, b_1 \not = b_2$, with $b_1, b_2 \in \wt{B}$}.}\\

\quad{\bf SUBCASE (1)} {\it Let $\Phi(b_1b_2) = \ol{b}\ol{a} \in B$}.

Let $M \in \Phi^{-1}(A)$ be a line incident with $b_1$. Let $c_1 \: \mathbf{I} \: M$, with $c_1 \in \wt{C}$ (we know that $M$ is incident with at least one point of $\wt{C}$). Assume, by way of contradiction, that $b_2$ is incident with a second line $M^\prime \in \Phi^{-1}(B)$. Then $M^\prime$ contains a point $a_1 \in \wt{A}$. Let $N^\prime$ be a line of $\Phi^{-1}(B)$ incident with $c_1$ and let $d_1 \in \wt{D}$ be incident with $N^\prime$. Then {\bf{d}}$(a_1, d_1)$ = 8, a contradiction. Hence $b_2$ is incident with a unique line of $\Phi^{-1}(B)$. Similarly $b_1b_2$ is the unique line of $\Phi^{-1}(B)$ incident with any point $b_i$ of $\wt{B}$ incident with $b_1b_2$. So $b_i$ is incident with $t$ lines of $\Phi^{-1}(A)$.

Assume, by way of contradiction, that $b_i \sim b^\prime, b^\prime  {\not \mathbf{I}}  \: b_1b_2, b^\prime \in \wt{B}$. Then $b_ib^\prime \in \Phi^{-1}(A)$. Let $b_i \sim c^\prime, b_i \not = c^\prime, b_1b_2 \not = b_ic^\prime \not = b_ib^\prime$. Then $b_ic^\prime \in \Phi^{-1}(A)$. We may assume that $c^\prime \in \wt{C}$. Let $c^\prime \sim d^\prime, c^\prime \not = d^\prime, {c^\prime d^\prime} \in \Phi^{-1}(B), d^\prime \in \wt{D}$. Further, let $a^\prime \sim b^\prime, a^\prime \not = b^\prime, b_ib^\prime \not = {a^\prime b^\prime}, {a^\prime b^\prime} \in \Phi^{-1}(B)$ and $a^\prime \in \wt{A}$. Then {\bf{d}}$(a^\prime, d^\prime)$ = 8, a contradiction. Hence each line $N \in \Phi^{-1}(A)$ incident with $b_i$ contains $s$ points of $\wt{C}$.

Let $b_1, b_2, \ldots, b_u$ be the points of $\wt{B}$ incident with $b_1b_2$. Then all points not incident with $b_1b_2$ but collinear with one of $b_1, b_2, \ldots, b_u$ belong to $\wt{C}$.

Assume that $u < s$, so that $b_1b_2$ is incident with at least two points $a_1, a_2$ of $\wt{A}$. Assume, by way of contradiction, that there is a second line $R$ of $\Phi^{-1}(B)$ incident with $a_1$. Let $r_1 \: \mathbf{I} \: R, r_1 \in \wt{B}, r_1 \sim r_2, r_1 \not = r_2, r_1r_2 \in \Phi^{-1}(A), r_2 \in \wt{C}$, and let $r_2 \sim d_1, r_2 \not = d_1, r_2d_1 \in \Phi^{-1}(B), d_1 \in \wt{D}$. Then {\bf{d}}$(a_2, d_1)$ = 8, a contradiction. So each point of $\wt{A}$ incident with $b_1b_2$ is incident with $t$ lines of $\Phi^{-1}(A)$. If there would be a point of $\wt{A}$ not incident with $b_1b_2$ but collinear with a point $a_i$ of $\wt{A}$ incident with $b_1b_2$, then it is easy to construct a point $d_2 \in \wt{D}$ with {\bf{d}}$(a_i, d_2)$ = 8, a contradiction. Hence each point not incident with $b_1b_2$ but collinear with some point of $\wt{A}$ incident with $b_1b_2$, belongs to $\wt{H}$.  

Assume again that $u < s$. Assume, by way of contradiction, that there exists a point $c^{\prime\prime} \in \wt{C}$ which is not collinear with a point of $\wt{B}$ incident with $b_1b_2$.
We consider two cases.

\begin{itemize}
\item[(1.1)] {\it {\bf{d}}$(c^{\prime\prime}, b_1b_2) \in \lbrace 3, 5 \rbrace$}.

Clearly {\bf{d}}$(c^{\prime\prime}, b_1b_2)$ = 3 is not possible (this follows also from the foregoing paragraph). So assume that {\bf{d}}$(c^{\prime\prime}, b_1b_2)$ = 5. First, let $c^{\prime\prime} \sim c^{\prime\prime\prime} \sim b_i, c^{\prime\prime\prime} \in \wt{C}, 1 \le i \le u$. Assume that $c^{\prime\prime}c^{\prime\prime\prime} \in \Phi^{-1}(A)$. Then $c^{\prime\prime\prime}b_i \in \Phi^{-1}(A)$. If $d^{\prime\prime} \in \wt{D}, c^{\prime\prime} \sim d^{\prime\prime}$, and $c^{\prime\prime}d^{\prime\prime} \in \Phi^{-1}(B)$, then {\bf{d}}$(a_1, d^{\prime\prime})$ = 8, a contradiction. Hence $c^{\prime\prime}c^{\prime\prime\prime} \in \Phi^{-1}(B)$ and $c^{\prime\prime\prime}b_i \in \Phi^{-1}(A)$. Interchange roles of $\lbrace b_1, b_2 \rbrace$ and $\lbrace c^{\prime\prime}, c^{\prime\prime\prime} \rbrace$, and of $\wt{B}$ and $\wt{C}$. Then by the first section of (1) we know that $c^{\prime\prime\prime}b_i$ is incident with $s$ points of $\wt{B}$, a contradiction as $c^{\prime\prime\prime}b_i$ is also incident with $s$ points of $\wt{C}$. Consequently $c^{\prime\prime} \sim w \sim a_j$, with $a_j \in \wt{A}$ incident with $b_1b_2$. Then by the foregoing paragraph $w \in \wt{H}$, which contradicts $w \sim c''$.

\item
[(1.2)] {\it {\bf{d}}$(c^{\prime\prime}, b_1b_2)$ = 7}.

First, let $b_i \sim c^{\prime\prime\prime} \sim l \sim c^{\prime\prime}, 1 \le i \le u$, with $c^{\prime\prime}, c^{\prime\prime\prime} \in \wt{C}$. By (1.1) we have $l \not \in \wt{C}$. So $l \in \wt{B} \cup \wt{D}$. First, let $l \in \wt{B}$. Let $b^{\prime\prime\prime} \: \mathbf{I} \: c^{\prime\prime\prime}l, l \not = b^{\prime\prime\prime} \not = c^{\prime\prime\prime}$.  Then by (1.1) we have $b^{\ppp} \in \wt{B}$. Let $a^{\pp} \sim b^{\ppp}, a^{\pp} \in \wt{A}$, and let $d^{\pp} \sim c^{\pp}, d^{\pp} \in \wt{D}$. Then {\bf{d}}$(a^{\pp} , d^{\pp})$ = 8, a contradiction. Hence $l \in \wt{D}$. Let $l^\prime \: \mathbf{I} \: lc^{\ppp}, l \not = l^\prime \not = c^{\ppp}$. By (1.1) $l^\prime \in \wt{D}$. Let $c^{\pp} \sim b^{\pppp} \sim a^{\ppp}, b^{\pppp} \in \wt{B}, a^{\ppp} \in \wt{A}$. Then {\bf{d}}$(l^\prime, a^{\ppp})$ = 8, a contradiction. Next, let $a_j \sim h^\prime \sim z \sim c^{\pp}$, with $a_j \: \mathbf{I} \: b_1b_2, a_j \in \wt{A}, h^\prime \in \wt{H}$. Then {\bf{d}}$(h^\prime, c^{\pp})$ = 4, clearly a contradiction.
\end{itemize} 

From (1.1) and (1.2) follows that for $u < s$ each point of $\wt{C}$ is collinear with some point of $\lbrace b_1, b_2, \ldots, b_u \rbrace$.

Now we will prove that for $u \le s$ each point of $\wt{A}$ is incident with the line $b_1b_2$. Assume, by way of contradiction, that $a^\ast \in \wt{A}$ is not incident with $b_1b_2$. Let $b_i \sim c_j \sim d_1, c_j \in \wt{C}, d_1 \in \wt{D}$. If the line $c_jd_1$ is incident with at least two points of $\wt{C}$, then, interchanging roles of $\wt{B}$ and $\wt{C}$, we see that $s$ points of $b_ic_j$ must belong to $\wt{B}$, a contradiction. Hence $c_jd_1$ is incident with $s$ points of $\wt{D}$. Let $d_2$ be a second point of $\wt{D}$ incident with $d_1c_j$. Then {\bf{d}}$(a^\ast, d_1)$ = {\bf{d}}$(a^\ast, d_2)$ = 6, recalling that $\mathbf{d}(a^*,d_1) = 8$, $\mathbf{d}(a^*,d_2) = 8$ are not possible, 
and so {\bf{d}}$(c_j, a^\ast)$ = 4. Let $c_{j^\prime} \: \mathbf{I} \: b_ic_j, c_{j^\prime} \in \wt{C}$. Then also {\bf{d}}$(c_{j^\prime}, a^\ast)$ = 4; this holds for the $s$ points $c_{j^\prime}$ of $\wt{C}$ incident with $b_ic_j$ as $\mathbf{d}(a^*,c_{j'}) \geq 4$. Hence $a^\ast \sim b_i$, and so $a^\ast \: \mathbf{I} \: b_1b_2$, a contradiction. Consequently every point of $\wt{A}$ is incident with the line $b_1b_2$. If $h_1 \in \wt{H}$, then any line of $\Phi^{-1}(A)$ incident with $h_1$ is incident with a point $a_j \in \wt{A}$. Hence each point of $\wt{H}$ is collinear with some point $a_j \: \mathbf{I} \: b_1b_2$.

\begin{quote}
{\footnotesize  
{\bf We summarize}. {\it Every point $b_i, i \in \lbrace 1, 2, \ldots, u \rbrace$, is incident with $t$ lines of $\Phi^{-1}(A)$ and just one line of $\Phi^{-1}(B)$. Every point not incident with $b_1b_2$ but collinear with some $b_i, i \in \lbrace 1, 2, \ldots, u \rbrace$, belongs to $\wt{C}$.} 
Assume now that $u < s$. Then every point $a_j \in \wt{A}$ incident with $b_1b_2$ is incident with $t$ lines of $\Phi^{-1}(A)$; each of these lines is incident with $s$ points of $\wt{H}$. Also, each point of $\wt{C}$ is collinear with some point of $\lbrace b_1, b_2, \ldots, b_u \rbrace$.
Finally, for $u \le s$ each point of $\wt{A}$ is incident with $b_1b_2$, and each point of $\wt{H}$ is collinear with some point of $\lbrace a_1, a_2, \ldots, a_{s + 1 - u} \rbrace$}. 
\end{quote} 

Assume that there is a point $b^\ast \! {\not\mathbf{I}} \: b_1b_2, b^\ast \in \wt{B}$. A line $V$ of $\Phi^{-1}(B)$ incident with $b^\ast$ contains a point $a_j \in \wt{A}$, with $a_j \: \mathbf{I} \: b_1b_2$; the other points incident with $V$ belong to $\wt{B}$.

\begin{itemize}
\item[(A)] {\it Assume that $\Phi^{-1}(\ol{a}\ol{b}) = \lbrace b_1b_2 \rbrace$}.

If $b^\ast \! {\not\mathbf{I}} \: b_1b_2, b^\ast \in \wt{B}$, then $b^\ast$ is incident with a line of $\Phi^{-1}(\ol{a}\ol{b})$, a contradiction. Hence each point of $\wt{B}$ is incident with $b_1b_2$. In this case $\wt{A}, \wt{B}$ define a partition of the set of all points incident with $b_1b_2, 1 < \vert \wt{B} \vert \le s, 1 \le \vert \wt{A} \vert \le s$. Also, $\wt{C}$ consists of all points not incident with $b_1b_2$ but collinear with a point of $\wt{B} = \lbrace b_1, b_2, \ldots, b_u \rbrace$, and $\wt{H}$ consists of all points not incident with $b_1b_2$ but collinear with a point of $\wt{A} = \lbrace a_1, a_2, \ldots, a_{s + 1 - u} \rbrace $. Now it easily follows that we have Case (a) in the statement of the theorem.

\item[(B)] {\it Assume that $\vert \Phi^{-1}(\ol{a}\ol{b}) \vert > 1$}.

First, let $\vert \wt{A} \vert > 1$ and let $U \in \Phi^{-1}(\ol{a}\ol{b}), U \not = b_1b_2$. So $U \in \Phi^{-1}(B)$. But by the summary,  we have that as $U$ is incident with a point of $\wt{A}$, the line $U$ belongs to $\Phi^{-1}(A)$, a contradiction.

Hence $\wt{A} = \lbrace a_1 \rbrace$. Then $\wt{B}$ consists of the points distinct from $a_1$ but incident with $r$ lines incident with $a_1, 1 \le r \le t$, and $\wt{H}$ consists of the points distinct from $a_1$ but incident with $t + 1 - r$ lines incident with $a_1$. Now it easily follows that we have Case (b) in the statement of the theorem.
\end{itemize}

\quad{\bf SUBCASE (2)} {\it Let $\Phi(b_1b_2) = \ol{b}\ol{c} \in A$}.

Let $b_1, b_2, \ldots, b_u  (u > 1)$ be the points of $\wt{B}$ incident with $b_1b_2$, and let $c_1, c_2, \ldots, c_{s + 1 - u}$ be the points of $\wt{C}$ incident with $b_1b_2$. By (1) we may assume that there is no point $b^\ast \in \wt{B}, b^\ast \not \in \lbrace b_1, b_2, \ldots, b_u \rbrace$, with $b_i \sim b^\ast$ and $b_ib^\ast \in \Phi^{-1}(B), i \in \lbrace 1, 2, \ldots, u \rbrace$. A line of $\Phi^{-1}(B)$ incident with $b_i, i \in \lbrace 1, 2, \ldots, u \rbrace$, contains $s$ points of $\wt{A}$.

Assume, by way of contradiction, that there is a second line $V \in \Phi^{-1}(A)$ incident with $b_i,  i \in \lbrace 1, 2, \ldots, u \rbrace$. Let $c^\ast \: \mathbf{I} \: V, c^\ast \in \wt{C}, d \sim c^\ast \sim b_i \sim b_j \sim a$ with $d \in \wt{D}, j \in \lbrace 1, 2, \ldots, u \rbrace, i \not = j, a \in \wt{A}$. Then {\bf{d}}$(a, d)$ = 8, a contradiction. So each $b_i, i \in \lbrace 1, 2, \ldots, u \rbrace$, is incident with $t$ lines of $\Phi^{-1}(B)$.

Let $c_j \sim d^\prime, d^\prime \in \wt{D}, j \in \lbrace 1, 2, \ldots, s + 1 - u \rbrace$. Assume that $c^\prime \in \wt{C}, c^\prime \: \mathbf{I} \: c_jd^\prime, c_j \not = c^\prime$. Interchanging roles of $\wt{B}$ and $\wt{C}$ in Case (a), then by Subcase (1) we have one of the cases in the statement of the theorem. So we may assume that the line $c_jd^\prime$ is incident with $s$ points of $\wt{D}$.

Assume that there is a point $d^{\pp} \in \wt{D}$ not collinear with one of the points $c_1, c_2, \ldots, c_{s + 1 - u}$. Let $a^\prime \sim a^{\pp} \sim b_1 \sim a^\prime$, with $a^\prime \not = a^{\pp}$ and $a^\prime, a^{\pp} \in \wt{A}$. Then {\bf{d}}$(a^\prime, d^{\pp})$ = {\bf{d}}$(a^{\pp}, d^{\pp})$ = 6. Then necessarily $d^{\pp}$ is collinear with some point $c_j (d^{\pp} \sim c_j \sim b_1)$ of $\lbrace c_1, c_2, \ldots, c_{s + 1 - u} \rbrace$, a contradiction. Hence $\wt{D}$ is the set of all points not incident with $b_1b_2$ but collinear with some point $c_j, j \in \lbrace 1, 2, \ldots, s + 1 - u \rbrace$.

Assume, by way of contradiction, that there is some point $c^{\pp} \in \wt{C}$ which is not incident with $b_1b_2$. Let $d^\ast \sim c^{\pp}$, with $d^\ast \in \wt{D}$. Then $d^\ast$ is collinear with some $c_j, j \in \lbrace 1, 2, \ldots, s + 1 - u \rbrace$. Interchanging again roles of $\wt{B}$ and $\wt{C}$, then by (a)(1) we may assume that the line $c^{\pp}d^\ast$ is incident with $s$ points of $\wt{D}$ (otherwise we have one of the two cases in the statement of the theorem). Each such point is collinear with some point of $\wt{C}$ incident with $b_1b_2$, a contradiction. So every point of $\wt{C}$ is incident with $b_1b_2$.

Now assume that $u < s$, that is, there are at least two points $c_1, c_2$ of $\wt{C}$ incident with $b_1b_2$. Assume, by way of contradiction, that the point $a^{\ppp} \in \wt{A}$ is not collinear with some point of $\lbrace b_1, b_2, \ldots, b_u \rbrace$. Let $c_1 \sim d^{\ppp} \sim d^{\pppp} \sim c_1, d^{\ppp} \not = d^{\pppp}$, and $d^{\ppp}, d^{\pppp} \in \wt{D}$. Then {\bf{d}}$(a^{\ppp}, d^{\ppp})$ = {\bf{d}}$(a^{\ppp}, d^{\pppp})$ = 6. So necessarily {\bf{d}}$(a^{\ppp}, c_1)$ = 4. Similarly, {\bf{d}}$(a^{\ppp}, c_2)$ = 4, clearly a contradiction. Hence every point of $\wt{A}$ is collinear with some point of $\lbrace b_1, b_2, \ldots, b_u \rbrace$.

Assume again that $u < s$. Interchanging roles of $\wt{B}$ and $\wt{C}$, the second paragraph of (a)(2) shows that each $c_j$ is incident with $t$ lines of $\Phi^{-1}(B), j = 1, 2, \ldots, s +1 -u$.

Let again be $u < s$. By interchanging roles of $\wt{B}$ and $\wt{C}$, we see, relying on a foregoing paragraph, that each point of $\wt{B}$ is incident with $b_1b_2$. Hence the set of all points incident with $b_1b_2$ is the set $\wt{B} \cup \wt{C}$.

Now it easily follows that for $u < s$ we have Case (a) in the statement of the theorem.

Consequently we now have to assume that $\vert \wt{C} \vert = 1$, so $\wt{C} = \lbrace c_1 \rbrace$. 

\begin{quote}
{\footnotesize 
{\bf We summarize what we know in this case}. {\it Each $b_i$ is incident with $t$ lines of $\Phi^{-1}(B), i = 1, 2, \ldots, u$; each such line contains $s$ points of $\wt{A}$. The unique point $c_1$ of $\wt{C}$ is incident with $t^\prime, 1 \le t^\prime \le t$, lines of $\Phi^{-1}(B)$; each such line contains $s$ points of $\wt{D}$. Each point of $\wt{D}$ is collinear with $c_1$.}} 
\end{quote} 

We proceed with the proof of the theorem. 

\begin{itemize}
\item[(A)] {\it Assume that $t^\prime = t$}.

Assume, by way of contradiction, that $a^\ast \in \wt{A}$ is not collinear with some point of $\lbrace b_1, b_2, \ldots, b_s \rbrace$. As {\bf{d}}$(a^\ast, d^\ast) \ne 8$ for each $d^\ast \in \wt{D}$, it follows that $a^\ast \sim b^\ast \sim c_1$ for some $b^\ast \in \wt{B}$ not incident with $b_1b_2$. Hence $b^\ast c_1 \in \Phi^{-1}(A)$, so $t^\prime < t$, a contradiction. Hence each point of $\wt{A}$ is collinear with some $b_i, i \in \lbrace 1, 2, \ldots, s \rbrace$.

Finally, we show that each point of $\wt{B}$ is incident with $b_1b_2$. Assume, by way of contradiction, that $b^{\ast\ast} \in \wt{B}, b^{\ast\ast}  \! {\not\mathbf{I}} \: b_1b_2$. Let $a^{\ast\ast} \sim b^{\ast\ast}, a^{\ast\ast} \in \wt{A}$. Then $a^{\ast\ast} \sim b_i, i \in \lbrace 1, 2, \ldots, s \rbrace$. By (a)(1) we may assume that $a^{\ast\ast}b^{\ast\ast}$ contains $s$ points of $\wt{A}$ (otherwise we have Case (a) or Case (b) in the statement of the theorem), and each of these points is collinear with some $b_j, j \in \lbrace 1, 2, \ldots, s \rbrace$, a contradiction. Hence each point of $\wt{B}$ is incident with $b_1b_2$.

As before it is now clear that we have Case (a) in the statement of the theorem.

\item[(B)] {\it Now assume $t^\prime < t$}.

Let $M \in \Phi^{-1}(A), M \not = b_1b_2$, be incident with $c_1$. Then $M$ is incident with $s$ points of $\wt{B}$. Also, each line of $\Phi^{-1}(B)$ containing a point of $\wt{B}$ is incident with $c_1$. Now it is clear that we have Case (b) in the statement of the theorem.
\end{itemize} 

\fbox{
\quad{\bf CASE (b)} {\it Assume that $c_1 \sim c_2, c_1 \not = c_2$, with $c_1, c_2 \in \wt{C}$}.}\\

Simllar to (a).\\

\fbox{
\quad{\bf CASE (c)} {\it No two points of $\wt{B}$ and no two points of $\wt{C}$ are collinear}.}\\

Let $b_1 \in \wt{B}, c_1 \in \wt{C}, b_1 \sim c_1$. Then necessarily $b_1, c_1$ are the only points incident with $b_1c_1$, so $s = 1$, a contradiction. \eop \\

%\medskip
%\subsection{Remark: epimorphisms involving locally finite octagons}
%\label{lf}

Call a thick generalized $n$-gon of order $(s,t)$ {\em locally finite} if exactly one of $s, t$ is finite. 

\begin{corollary}
\label{lf}
Let $\gamma: \Gamma \mapsto \Delta$ be an epimorphism from the thick locally finite generalized octagon $\Gamma$ of order $(s,t)$ with $t$ finite, to the thin octagon $\Delta$ of order $(s',1)$. Then $s' = 1$ and we have the conclusion of { Theorem \ref{JATGO}}.  
\end{corollary}

{\em Proof}.\quad
The part of the proof of Theorem \ref{JATGO} that shows that $s' = 1$ only uses the fact that the number of lines incident with a point of $\Gamma$ is finite. \eop \\

 \section{Locally finitely chained and generated generalized polygons}
  \label{chain}

%Recall that a group $G$ is called {\em locally finite} if every finite subset of $G$ generates a finite subgroup. 

Let $S$ be any point set of a generalized $n$-gon $\Gamma$; then as in \cite{part3}, $\langle S \rangle$ by definition is the intersection of all (thin and thick) sub $n$-gons that contain $S$. We call $\langle S \rangle $ the subgeometry {\em generated by $S$}. Note that this not necessarily is a (thin or thick) generalized $n$-gon itself, but generically it is. 
 
Call a thick generalized $n$-gon {\em locally finitely generated} if the following property holds: 
\begin{quote}
for every finite point subset $S$ which generates a (possibly thin) sub $n$-gon $\langle S \rangle$, we have that $\langle S \rangle$ is finite.  
\end{quote}
Finite generalized $n$-gons are trivially locally finitely generated. \\

Call a generalized $n$-gon $\Gamma$ {\em locally finitely chained} if there exists a chain of finite point subsets
\[  S_0 \subseteq S_1 \subseteq \ldots \subseteq S_i \subseteq \ldots                    \]       
indexed over the positive integers, such that each $S_j$ generates a finite  (possibly thin) sub $n$-gon, and such that 
\[  \bigcup_{i \geq 0}\Big\langle S_i \Big\rangle\ =\ \Gamma.  \]
If $\Gamma$ is not finite, it is easy to see that we may assume that the chain 
\[  \langle S_0 \rangle \subseteq \langle S_1\rangle \subseteq \ldots \subseteq \langle S_i \rangle \subseteq \ldots                    \]     
is strict.\\

A very simple key observation is the following lemma taken from \cite{part3}, in which we suppose that  the considered locally finitely generated and locally finitely 
chained $n$-gons are not finite (to avoid trivialities). 

\begin{lemma}[Thas and Thas \cite{part3}]
\label{lem1}
\begin{itemize}
\item[{\rm (a)}]
Every thick locally finitely generated generalized $n$-gon contains thick locally finitely chained sub $n$-gons. 
\item[{\rm (b)}]
Every thick finitely chained generalized $n$-gon is locally finitely generated.
\item[{\rm (c)}]
A thick generalized $n$-gon is locally finitely generated if and only if its point-line dual is.
\item[{\rm (d)}]
A thick generalized $n$-gon is locally finitely chained if and only if its point-line dual is.
\end{itemize}
\end{lemma}

Note that $\cup_{i \geq 0}\Big\langle S_i \Big\rangle$ is countably infinite or finite.\\

\begin{observation}[Thas and Thas \cite{part3}]
\begin{itemize}
\item[{\rm (a)}]
Let $\Gamma$ be a thick locally finitely chained generalized $n$-gon. Then $n \in \{ 3, 4, 6, 8 \}$.
\item[{\rm (b)}]
Let $\Gamma$ be a thick locally finitely generated generalized $n$-gon. Then $n \in \{ 3, 4, 6, 8 \}$.
\end{itemize} 
\end{observation}

(Due to their very definition, both types of polygon share many properties with finite polygons.) 

Before proceeding, we introduce the following two properties for a point-line incidence structure $\Gamma = (\mP,\mL,\I)$:
\begin{itemize}
\item[(Ia)]
every two elements of $\mP \cup \mL$ can be joined  by at most one path of length $< n$; 
\item[(Ib)]
every two elements of $\mP \cup \mL$ can be joined by at least one path of length $\leq n$. 
\end{itemize}

Call a point-line incidence geometry $\Gamma = (\mP,\mL,\I)$ {\em firm} if each element is incident with at least two different elements. 
Then by Van Maldeghem \cite{POL2}, $\Gamma$ is a weak generalized $n$-gon if and only if it is firm, and both (Ia) and (Ib) are satisfied. 

\begin{theorem}
Let $\Gamma$ be a thick locally finitely chained generalized $8$-gon, and let $\gamma: \Gamma \mapsto \Delta$ be an epimorphism, with $\Delta$ a thin $8$-gon of order $(s,1)$. Then $s = 1$ and the conclusions of Theorem \ref{JATGO} hold.
\end{theorem}

{\em Proof.}\quad 
Assume, by way of contradiction, that $s \ne 1$. 
Now define $m$ such that $\gamma(\Big\langle S_m \Big\rangle)$ strictly contains an ordinary $8$-gon $\Delta'$ as subgeometry, and such that $\Big\langle S_m \Big\rangle$ is finite and thick. Since  $\gamma(\Big\langle S_m \Big\rangle) =: \Delta''$ is contained in $\Delta$, we have that (Ia) is automatically satisfied in $\Delta''$. And (Ib) can be verified in $\Delta''$ by applying $\gamma^{-1}$. Since $\Delta''$ contains an ordinary $8$-gon, it now follows that it is firm, so it is a weak generalized $8$-gon. Also, 
\[ \gamma:\ \Big\langle S_m \Big\rangle\ \mapsto\ \Delta''              \]
is an epimorphism. By \cite[Theorem 3.1]{POL2} (see also \cite[Theorem 1.6.2]{POL}), $\Delta''$ is either 
\begin{itemize}
\item
the point-line dual of the double of a generalized quadrangle, or 
\item
the point-line dual of a degenerate octagon $\bO$ which consists of opposite points $a$ and $b$, joined by $r \geq 2$ paths of length $8$, or 
\item
the point-line dual of the quadruple of a digon (in which case the minimal distance between thick elements in $\Delta''$ is $4$). 
\end{itemize}

The point-line dual of $\bO$ has at most two thick lines, so we can take $m$ large enough so that $\gamma(\langle S_m \rangle)$ contains more than two thick lines. As $s > 1$ by assumption, we can also take $m$ large enough such that  $\gamma(\Big\langle S_m \Big\rangle)$ contains thick elements at distance $2$. 

By the proof of Theorem \ref{JATGO}, the statement now follows (as $\Delta''$ is the point-line dual of the double of a generalized quadrangle).
 \eop \\
 
% \textcolor{blue}{HIER IS HET ESSENITEEL IN THE CASE $n = 6$ DAT $\Delta$ EEN DUBBELE IS VAN EEN PV. DUS DAT MOETEN WE WEL ZEKER ZIJN.}

\begin{corollary}
Let $\Gamma$ be a thick locally finitely generated generalized $8$-gon, and let $\gamma: \Gamma \mapsto \Delta$ be an epimorphism, with $\Delta$ a thin $8$-gon of order $(s,1)$. Then $s = 1$ and the conclusions of Theorem \ref{JATGO} hold. \eop \\
\end{corollary}

\subsection*{EXAMPLES}

The property of being locally finitely chained / generated appears to be rather strong, but there are surprisingly many interesting examples of such polygons with much structure as 
we have seen in \cite{part3}. In this section, we describe examples using finite classical octagons. 

Suppose $\F_q$ is a finite field, and let $\overline{\F_q}$ be an algebraic closure of $\F_q$. Let $q = p^m$ for the prime $p$. Then for each positive integer $n$, $\overline{\F_q}$ has a unique subfield isomorphic to $\F_{p^n}$ (we denote it also as $\F_{p^n}$).  Then 
\[ \bigcup_{i \geq 1}\F_{p^i}\ =\ \overline{\F_q}.    \]

In particular, $\F_{p^n}$ is a subfield of $\F_{p^m}$ if and only if $n \vert m$, and there is only one subfield of $\F_{p^m}$ which is isomorphic to $\F_{p^n}$. For each field 
$\F_{p^n}$, the map $\gamma_p: \F_{p^n} \mapsto \F_{p^n}: x \mapsto x^p$ is a field automorphism which is called the {\em Frobenius automorphism}. 

Consider an infinite set of pairs $\Big\{ (\F_{2^{2n_i + 1}}, \sigma_i )    \Big\}_i$, where $\sigma_i$ is an automorphism of  $\F_{2^{2n_i + 1}}$ such that $\sigma_i^2(x) = x^2$ for all 
$x$ in the field (that is, $\sigma_i^2$ is the Frobenius automorphism), and so that if $i < j$, then $2n_i + 1$ divides $2n_j + 1$. In that case, $\F_{2^{2n_i + 1}}$ is a subfield of $\F_{2^{2n_j + 1}}$ and $\sigma_j$ fixes $\F_{2^{2n_i + 1}}$, so it induces an automorphism of the latter. Note that $\sigma_i$ is defined as follows: $\sigma_i: \F_{2^{2n_i + 1}}\mapsto \F_{2^{2n_i + 1}}: x \mapsto x^{n_i + 1}$, and also note that the actions of $\sigma_i$ and $\sigma_j$ coincide on $\F_{2^{2n_i + 1}}$ so that the system $\Big\{ (\F_{2^{2n_i + 1}}, \sigma_i )    \Big\}_i$ is ``compatible.'' 

For each pair $(k_i := \F_{2^{2n_i + 1}}, \sigma_i )$ we construct a Ree-Tits octagon $\bO(k_i,\sigma_i)$ as defined in \cite[section 2.5]{POL}. 

We will use the following theorem:

\begin{theorem}[Joswig and Van Maldeghem \cite{JosHVM}]
\label{JHVM}
Let $\bO(k,\sigma)$ be a Ree-Tits octagon. If $k'$ is a subfield of $k$ for which ${k'}^\sigma \subseteq k'$, then $\bO(k',\sigma_{\vert k'})$ naturally defines a sub Ree-Tits octagons by ``restricting coordinates.''  Vice versa, if $\bO'$ is a proper thick suboctagon of $\bO(k,\sigma)$, then $\bO'$ is a Ree-Tits octagon which arises in this way (by field reduction). 
\end{theorem}

By Theorem \ref{JHVM} we know that the set  $\Big\{ (k_i, \sigma_i )    \Big\}_i$ actually defines a chain of finite Ree-Tits octagons, and obviously the union of these octagons 
is again a (non-finite) Ree-Tits octagon $\bO(\widehat{k},\widehat{\sigma})$, where $\widehat{k}$ is the union of the fields $k_i$, and $\widehat{\sigma}$ is the automorphism of 
$\widehat{k}$ defined by the local actions of the $\sigma_i$ induced on the subfields $k_i$.  (Note that if $a \in \widehat{k}$, then $a$ is contained in some subfield $k_m$, 
and $\widehat{\sigma}(a) = \sigma_m(a)$, so that trivially $\widehat{\sigma}^2$ is the Frobenius automorphism of $\widehat{k}$.) Each of $\bO(\widehat{k},\widehat{\sigma})$, $\widehat{k}$ and $\widehat{\sigma}$ can be seen as direct limits over the system $\Big\{ 2n_i + 1 \Big\}_i$.

It is very easy to find systems $\Big\{ 2n_i + 1 \Big\}_i$ as above. We describe a couple of them, but the possibilities are obviously numerous.
\begin{itemize}
\item
Let $N$ be any odd positive integer different from $1$, and let $\Big\{ m_i \Big\}_{i}$ be any strictly ascending chain of strictly positive integers; then define $2n_i + 1$ as  
$N^{m_i}$. 
\item
Enumerate the odd primes as $p_1, p_2, p_3, \ldots$. Then define $2n_1 + 1 = p_1$, $2n_2 + 1 = p_1p_2$, $2n_3 + 1 = p_1p_2p_3$, etc.
\item
Variations of the above. 
\end{itemize}

Note that it is very easy to construct nonisomorphic fields $\widehat{k}$ (and hence nonisomorphic octagons) using these constructions. For instance, let $N$ be any prime number; then 
the field $\widehat{k}_1$ constructed in the first example contains subfields isomorphic to $\F_{2^{N^m}}$ with $m > 1$, while the field $\widehat{k}_2$ constructed in the second example does not. 

All the octagons $\bO(\widehat{k},\widehat{\sigma})$ are locally finitely chained, and as they are countably infinite, they are also locally finitely generated. 

\begin{remark}{\rm 
Note that the only known thick finite generalized octagons are Ree-Tits octagons. If there would be no others (as has been conjectured by some), then the examples constructed above essentially describe all (infinite) 
locally finitely chained generalized octagons. For, let $\bO$ be a locally finitely chained thick generalized octagon. Then there is a chain 
\[  S_0 \subseteq S_1 \subseteq \ldots \subseteq S_i \subseteq \ldots                    \]       
indexed over the positive integers, such that each $S_j$ generates a finite (possibly thin) sub $8$-gon, and such that 
\[  \bigcup_{i \geq 0}\Big\langle S_i \Big\rangle\ =\ \bO.  \]
It is obvious that from some index $m$ on, $\Big\langle S_n \Big\rangle$ will be thick if $n \geq m$. So we can as well use only thick sub octagons (with the same notation as above). 
Since all of these are finite, we conjecturally end up with a chain of Ree-Tits octagons $\Big\langle S_i \Big\rangle = \bO(k_i,\sigma_i)$, where $u < v$ (positive integers) 
implies that $k_u$ is a subfield of $k_v$ and $\sigma_v$ induces $\sigma_u$ in $k_u$, and where $k_i \cong \F_{2^{2n_i + 1}}$ for some positive integer $n_i$. 
It follows that $\bO$ is of the type described above.   
}
\end{remark}

%Here, $\varinjlim(\cdot)$ denotes direct limit. 

%\begin{remark}[Planes over $\overline{\F_p}$]
%\end{remark}

%\subsection*{Free locally finitely chained polygons}

%\textcolor{blue}{DIT MOET IK NOG TOEVOEGEN.}

\medskip
\section{Counter examples in the infinite case}
\label{counter}

In this section we show that Theorem \ref{JATGO} does not hold for infinite generalized octagons, by constructing epimorphisms from 
thick infinite generalized $n$-gons to generalized $n$-gons with parameters $(s',1)$, $s' > 1$. Here, $n \geq 3$. \\

We first need some more terminology. Suppose that $\Upsilon$ is a point-line geometry; we measure distances using the incidence graph $\mathcal{I}$ of $\Upsilon$. Call a path $x = x_0, x_1,\ldots,x_n = y$ of length $n$ between vertices $x$ and $y$ {\em nonstammering} if $x_{i - 1} \ne x_{i + 1}$ for all $i \in \{1,2,\ldots,n - 1\}$. A {\em circuit} is a finite nonstammering closed path. The {\em girth} of $\Upsilon$ is the length of a minimal circuit (length $0$ is not allowed in this definition). It is $\infty$ by definition if no nontrivial circuits exist.  If the girth $n$ of $\Upsilon$ is finite, it is easy to see that it must be even, and we call $n/2$ the {\em gonality} of $\Upsilon$.

\begin{proposition}
Let $\Gamma$ be any generalized $n$-gon of order $(s',1)$ with $s' > 1$ finite or countable, and $n \geq 3$. Then there exist epimorphisms ${\varepsilon}:\ \mA\ \mapsto\ \Gamma$, with $\mA$ a thick generalized $n$-gon.  
\end{proposition}

{\em Proof.}\quad 
Suppose that $A = A_0$ is a connected point-line geometry of finite gonality $m \geq n$ which is not a weak generalized $n$-gon, and where $n \geq 3$,  and let $\epsilon: A \mapsto \Gamma$ be an 
epimorphism (of point-line incidence geometries), where $\Gamma$ is a generalized $n$-gon of order $(s',1)$, where $s' > 1$ is arbitrary but finite or countable. We also suppose that $A$ has a finite or countable number of points and/or lines. (The fact that $A$ is connected is not important: starting from a disconnected geometry, one can make it 
connected by adding appropriate paths between the connected components. )

We will freely construct a generalized $n$-gon $\overline{A}$ over $A$ and an epimorphism $\overline{\epsilon}: \overline{A} \mapsto \Gamma$. 

Suppose $A_i$ and the epimorphism $\epsilon_i: A_i \mapsto \Gamma$ were constructed in a previous step, and define $A_{i + 1}$ as follows: if $x$ and $y$ are elements at distance $n - 1$ in $A_i$, then we add 
a completely new path $\gamma(x,y)$ of length $n + 1$ between $x$ and $y$. Now define $\epsilon_{i + 1}$ as follows:
\begin{itemize}
\item
if the distance between $\epsilon_i(x)$ and $\epsilon_i(y)$ is smaller than $n - 1$ in $\Gamma$, then map $\gamma(x,y)$ surjectively on the unique shortest path between $\epsilon_i(x)$ and $\epsilon_i(y)$;
\item
if the distance between $\epsilon_i(x)$ and $\epsilon_i(y)$ equals $n - 1$ in $\Gamma$, then map $\gamma(x,y)$ surjectively on an arbitrary path of length $n + 1$ between $\epsilon_i(x)$ and $\epsilon_i(y)$.
\end{itemize}

Now let $\overline{A} = \cup_{n \geq 1}A_n$ and $\overline{\epsilon} = \cup_{n \geq 1}\epsilon_n$ (where we identify a function with its graph). Then 
\[  \overline{\epsilon}:\ \overline{A}\ \mapsto\ \Gamma      \]
obviously is an epimorphism. 

By \cite[1.3.122]{POL}, it follows that $\overline{A}$ is a weak generalized $n$-gon; if $A$ already contained an ordinary $(n + 1)$-gon, it follows that $\overline{A}$ would be a (thick) generalized $n$-gon. In such case, setting $\Big(\mA,\varepsilon\Big) = \Big(\overline{A},\overline{\epsilon} \Big)$ yields the desired examples. 
\eop \\

Note that the construction above also works for thick targets $\Gamma$. 

Note also that it is easy to construct legitimate  begin-configurations $(A,\epsilon)$ (we will leave this exercise to the interested reader). \\

\begin{remark}{\rm 
In Gramlich and Van Maldeghem \cite{POL2}, it is proved that starting from any (thick) generalized $n$-gon $\Gamma$ ($n \geq 2$), one can construct a free generalized $n$-gon $\overline{A}$ and an epimorphism 
\[  \overline{\epsilon}:\ \overline{A}\ \mapsto\ \Gamma.      \]

The begin configuration in \cite{POL2} is different, and the target is thick, 
but of course the essence is the same. }
\end{remark}

%CHECK NOTES JAT: can the finiteness arguments be taken over? (In this case, we suppose there is a finite number of lines on any point, since we project onto a grid.) 

%TONGUE-in-CHEEK: can we construct such epimorphisms in case of hypothetical examples ? 

\end{document}